\documentstyle[amsfonts,12pt]{article}
\newcommand{\lon}{\longrightarrow}
\newcommand{\rar}{\rightarrow}

\newcommand{\Z}{{\Bbb Z}}
\newcommand{\p}{{\partial}}

\newcommand{\ot}{\otimes}
\newcommand{\Id}{\mbox{Id}}
\newcommand{\Beq}{\begin{equation}}
\newcommand{\Eeq}{\end{equation}}
\newcommand{\Beqr}{\begin{eqnarray}}
\newcommand{\Eeqr}{\end{eqnarray}}
\newcommand{\Beqrn}{\begin{eqnarray*}}
\newcommand{\Eeqrn}{\end{eqnarray*}}
\newcommand{\Ba}{\begin{array}}
\newcommand{\Ea}{\end{array}}
\newcommand{\Barr}{\begin{array}}
\newcommand{\Earr}{\end{array}}
\newcommand{\Bi}{\begin{itemize}}
\newcommand{\Ei}{\end{itemize}}
\newcommand{\Bc}{\begin{center}}
\newcommand{\Ec}{\end{center}}
\newcommand{\al}{\alpha}
\newcommand{\be}{\beta}

\newcommand{\la}{\lambda}

\newcommand{\tv}{\tilde{v}}
\newcommand{\tl}{\tilde{l}}
\newcommand{\tk}{\tilde{k}}

\newcommand{\tj}{\tilde{j}}

\newcommand{\bp}{\bar{\partial}}
\newcommand{\bbj}{\bar{j}}

\newcommand{\bbs}{\bar{s}}
\newcommand{\bbt}{\bar{t}}
\newcommand{\bbz}{\bar{z}}
\newcommand{\Ker}{\mbox{Ker}}
\newcommand{\vst}{\vspace{2 mm}}
\newcommand{\vsf}{\vspace{4 mm}}
\newcommand{\vse}{\vspace{8 mm}}
\newtheorem{theorem}{Theorem}[section]

\newtheorem{lemma}[theorem]{Lemma}

\newtheorem{assumption}[theorem]{Assumption}

\begin{document}

\title{ Strong homotopy algebras \\ of
a K\"ahler manifold}

\author{ S.A.\ Merkulov}
\date{}
\maketitle

\begin{abstract}
It is shown that any compact K\"ahler manifold $M$ gives canonically
rise to two strong homotopy  algebras, the first one being associated with the
Hodge theory of the de Rham complex and the second one  with the
Hodge theory of
the Dolbeault complex. In these algebras the product
of two
 harmonic differential forms is again harmonic.

\vst

If $M$ happens to be a Calabi-Yau manifold, there exists a third
strong homotopy  algebra closely related to the Barannikov-Kontsevich
extended moduli space of complex structures.

\end{abstract}

\vse

\section{Introduction}
Strong homotopy algebras have been introduced by Stasheff \cite{Stasheff1}
more than 30 years ago in the context of topological $H$-spaces.
Since that time this remarkable structure (and its
cousins)
has been found at a number of unexpected places, for example
in the string field theory, in the topological conformal field theory,
in the Morse theory and
in the symplectic Floer theory (see, e.g.,
\cite{Stasheff2,Stasheff3,Fu,Ko1,Ko2} and references
cited therein).

\vst

In this paper we show that in the heart of the classical Hodge theory of
K\"ahler manifolds (see, e.g., \cite{GH})
there lies a strong homotopy algebra which
 ``explains''
the following  imperfect
behavior of $d^*$-closed differential forms: the wedge
product, $\al\wedge \be$, of $d^*$-closed forms $\al$ and $\be$ is {\em not}\,
$d^*$-closed in general. This, of course, can be easily fixed by defining a new
product $\circ$,
$$
\al\circ \be:= [\al\wedge \be]_{{\mathrm Ker}d^*},
$$
$[\  ]_{{\mathrm Ker}d^*}$ being the natural projection to the 
${\mathrm Ker} d^*$ constituent,
but at a price --- the product $\circ$ fails to be associative. Among the
 main results of the paper is an observation that $\circ$ does satisfy the
{\em higher order associativity conditions}\, thereby unveiling the structure
of strong homotopy  algebra in the
Hodge theory of K\"ahler manifolds. The point is that we are able 
to write down {\em explicitly}\, all the higher order products
in terms of the wedge product, Green function and the K\"ahler form.

\vst

In fact there are at least two strong homotopy algebras associated with
a compact K\"ahler manifold $M$. The first one is real and is associated
with the Hodge theory of the de Rham complex $(\Omega^{\bullet}M, d)$;
the second one is complex (and is not complexification of the first) and
is associated with the Hodge theory  of
the Dolbeault complex $(\Omega^{\bullet,\bullet}M, \bar{\p})$.

\vst

If $M$ happens to be a Calabi-Yau manifold, there exists a third
strong homotopy algebra closely related to the Barannikov-Kontsevich
extended moduli space of complex structures \cite{BaKo}.

\vst

The paper is organized as follows. In Sect.\ 2 we recall the definition
of a strong homotopy algebra. In Sect.\ 3 we consider specific
algebraic data  which gives rise to a particular strong
homotopy algebra;
a remarkable feature of the construction is that it provides {\em explicit}\,
formulae for all the higher homotopies. In Sect.\ 4 we consider a
compact K\"ahler manifold and prove the results mentioned above.

\vse

%%%%%%%%%%%%%%%%%%%%%%%%%%%%%%%%%%%%%%%%%%%%%%%%%%%%%%%%
\section{ Strong homotopy algebras}

A {\em strong homotopy algebra}, or shortly $A_{\infty}$-{\em algebra},
is by definition a vector superspace $V$ equipped with linear maps,
$$
\Ba{rcccc}
\mu_k: & \ot^k V & \lon & V & \\
& v_1\ot\ldots \ot v_k & \lon & \mu_k(v_1, \ldots, v_k), & \ \ \ \ k\geq 1,
\Ea
$$
of parity $\tilde{k}:=k\bmod 2\Z$ satisfying, for any $n\geq 1$ and
 any $v_1, \ldots, v_n \in V$,
 the following
{\em higher order associativity conditions},
\Beq \label{id}
\sum_{k+l=n+1}\sum_{j=0}^{ k-1} (-1)^r \mu_k\left(v_1,\ldots,v_j,
\mu_l(v_{j+1},\ldots, v_{j+l}), v_{j+l+1}, \ldots, v_n\right)=0,
\Eeq
where $r=\tilde{l}(\tv_1 +\ldots + \tv_j) + \tj(\tl-1) + (\tk-1)\tl$
and $\tv$ denotes the parity of $v\in V$.

\vse

Denoting $dv_1:=\mu_1(v_1)$ and $v_1\circ v_2 := \mu_2(v_1,v_2)$, one may
depict explicitly
the first  four floors of the infinite tower of higher order
associativity conditions   as follows
\begin{description}
\item[$n=1$:]\ \  $d^2=0$,\\
\item[$n=2$:]\ \  $d(v_1\circ v_2)= (dv_1)\circ v_2 + (-1)^{\tv_1} v_1\circ
(dv_2)$,\\
\item[$n=3$:]\ \ $v_1\circ (v_2\circ v_3) - (v_1\circ v_2)\circ v_3 =
d \mu_3(v_1,v_2,v_3) + \mu_3(dv_1, v_2, v_3) + (-1)^{\tv_1} \mu_3(v_1, dv_2,
v_3) + \mbox{\hspace{5mm}$(-1)^{\tv_1+ \tv_2}\mu_3(v_1,v_2,dv_3)$}$,\\
\item[$n=4$:] \ \ $\mu_3(v_1,v_2,v_3)\circ v_4 - \mu_3(v_1\circ v_2, v_3, v_4)
+ \mu_3(v_1,v_2\circ v_3,v_4) - \mu_3(v_1,v_2,v_3\circ v_4)
+ (-1)^{\tv_1}v_1\circ \mbox{\hspace{5mm}$\mu_3(v_2,v_3,v_4)$} =
d\mu_4(v_1,v_2,v_3,v_4) - \mu_4(dv_1,v_2,v_3,v_4) - (-1)^{\tv_1}
\mu_4(v_1,dv_2,v_3,v_4) - \mbox{\hspace{4mm}
$(-1)^{\tv_1+\tv_2}\mu_4(v_1,v_2, dv_3,v_4)$}
- (-1)^{\tv_1+\tv_2 + \tv_3}\mu_4(v_1,v_2, v_2, v_3,dv_4)$.
\end{description}

Therefore $A_{\infty}$-algebras with $\mu_k=0$ for $k\geq 3$
are nothing but the differential associative superalgebras with the
differential $\mu_1$ and the associative multiplication defined by
$\mu_2$. If, furthermore, $\mu_1=0$, one recovers  the usual associative
superalgebras.

\vst

The notion of $A_{\infty}$-algebra is a very natural
extension of the usual concept of associative superalgebra.
The following well-known fact \cite{PS} may serve as a  confirmation of
this  statement:
the moduli space of infinitesimal deformations of an associative
superalgebra $A$ within the class of associative superalgebras is
isomorphic to the second
 Hochschild cohomology group $\mbox{\rm Hoch}^2(A,A)$, while
the moduli space of infinitesimal deformations of the same $A$ within the class
of $A_{\infty}$-algebras can be identified with the full
Hochschild cohomology  $\mbox{\rm Hoch}^*(A,A)$. In this sense
$A_{\infty}$-algebra is a ``final'' concept.

\vse

%%%%%%%%%%%%%%%%%%%%%%%%%%%%%%%%%%%%%%%%%%%%%%%%%%%%%%%%
\section{ An explicit construction of an $A_{\infty}$-algebra}

Let $(V,d)$ be a differential associative superalgebra,
with $d$ denoting the differential, i.e.\ an odd linear map $d: V\lon V$
satisfying the Leibnitz identity
$$
d(v_1 \cdot v_2)= (dv_1)\cdot v_2 + (-1)^{\tv_1}v_1\cdot (dv_2),
$$
for any $v_1,v_2\in V$.

\vst

Let  $W$ be a sub complex of $(V,d)$, i.e.\ a vector subspace $W\subset V$
invariant under $d$. Note that we do {\em not}\, assume that
$W$ is a subalgebra of $V$. Instead we make the following
\begin{assumption}\label{assumption}
There exists an odd operator
$$
Q: V \lon V
$$
such that for any $v\in V$ the element $(1-[d,Q])v$ lies in the subspace
$W$, where $[\, , \, ]$ is the supercommutator.
\end{assumption}

The resulting datum $(W\subset V,d,Q,\, \cdot \, )$ is almost the same
as the one considered by Gugenheim and Stasheff in \cite{GS}
except that we do not assume that the operator
$$
P\equiv\Id-[d,Q]:\, V \lon W
$$
is identity when restricted to $W$; moreover, it may
not be even a surjection.

\vst

Since $[d,P]=0$, the product
$$
v_1\circ v_2 := (1-[d,Q])(v_1\cdot v_2)
$$
makes $(W,d,\circ)$ into a differential superalgebra which
is {\em not}\, associative in general. Under the additional
assumption that $P|_{W}=\Id$ it was proved in \cite{GS}
that $\circ$ does satisfy the higher associativity conditions.
We shall give a new proof of this fact (under weaker
assumption~\ref{assumption}); moreover, we shall be able
to compute all the higher homotopies $\mu_k$ {\em explicitly}\,
in terms of $d$,  $Q$ and $\cdot$ only.

\vst

First we  define a series of linear maps
$$
\la_n: \ot^n V \lon V, \ \ \ \ \ n\geq 2,
$$
starting with
$$
\la_2(v_1,v_2) := v_1\cdot v_2
$$
and then recursively, for $n\geq 3$,
\Beqr
\la_n(v_1,\ldots,v_n) & := & (-1)^{n-1} [Q\la_{n-1}(v_1,\ldots,v_{n-1})]
\cdot v_n -
(-1)^{n\tv_1}v_1 \cdot [Q\la_{n-1}(v_2,\ldots, v_n)] \nonumber\\
&& -\,\sum_{ k+l=n+1\atop k,l\geq 2 }(-1)^{k+(l-1)(\tv_1+\ldots+
\tv_k)} [Q\la_k(v_1, \ldots, v_k)] \cdot [Q\la_l(v_{k+1},\ldots,v_n)].
\nonumber \\
&& \label{la1}
\Eeqr
For example,
\Beqrn
\la_3(v_1,v_2,v_3) &=& [Q\la_2(v_1,v_2)]\cdot v_3 - (-1)^{\tv_1}
v_1\cdot [Q\la_2(v_2,v_3)], \\
\la_4(v_1,v_2,v_3,v_4) &=& -[Q\la_3(v_1,v_2,v_3)]\cdot v_4 -
(-1)^{\tv_1+\tv_2}[Q\la_2(v_1,v_2)]\cdot[Q\la_2(v_3,v_4)] \\
&& - v_1\cdot[Q\la_3(v_2,v_3,v_4)]\\
\la_5(v_1,v_2,v_3,v_4,v_5) &=& [Q\la_4(v_1,v_2,v_3,v_4)]\cdot v_5
+ (-1)^{\tv_1 +\tv_2 +\tv_3}[Q\la_3(v_1,v_2,v_3)]\cdot [Q\la_2(v_4,v_5)]\\
&& - [Q\la_2(v_1,v_2)]\cdot[Q\la_3(v_3,v_4,v_5)]
- (-1)^{\tv_1} v_1\cdot [Q\la_4(v_2,v_3,v_4,v_5)].
\Eeqrn

\vst

Setting formally
$$
\la_1:= -Q^{-1}
$$
(which makes sense because in all our formulae below $\la_1$ enters
in the combination $Q\la_1=-\Id$), one can rewrite the recursion (\ref{la1})
as follows
\Beq\label{la2}
\la_n(v_1,\ldots,v_n) = -\sum_{ k+l=n+1\atop k,l\geq 1 }(-1)^{k+(l-1)(\tv_1+\ldots+
\tv_k)} [Q\la_k(v_1, \ldots, v_k)] \cdot [Q\la_l(v_{k+1},\ldots,v_n)],
\Eeq
where now $n\geq 2$.

\begin{lemma}\label{lemma1}
The tensors $\la_k$, $k\geq 2$, satisfy the identities
$$
\Phi_n(v_1,\ldots,v_n) \equiv \sum_{k+l=n+1\atop k,l\geq 2}\sum_{j=0}^{k-1}
(-1)^r \la_k\left(v_1,\ldots,v_j, \la_l(v_{j+1},\ldots,v_{j+l}), v_{j+l+1},
\ldots, v_n\right)=0,
$$
$$
r=l(\tv_1+\ldots \tv_j) + j(l-1) + (k-1)l,
$$
for any $n\geq 3$ and any $v_1,\ldots,v_n\in V$.
\end{lemma}

\vst

\noindent{\bf Proof}. First we split
$\Phi_n$,
\Beqr
\Phi_n(v_1,\ldots,v_n) &=& \sum_{k+l=n+1\atop k,l\geq 2}(-1)^{(k-1)l}
\la_k\left(\la_l(v_{1},\ldots,v_{l}), v_{l+1},
\ldots, v_n\right) \, + \nonumber \\
&& \sum_{k+l=n+1\atop k,l\geq 2}
(-1)^{l(\tv_1+\ldots+\tv_{k-1}) + k-1} \la_k\left(v_1,\ldots,v_{k-1},
\la_l(v_{k},\ldots, v_n)\right) \, + \nonumber \\
&& \sum_{k+l=n+1\atop k,l\geq 2}\sum_{j=1}^{k-2}
(-1)^r \la_k\left(v_1,\ldots,v_j, \la_l(v_{j+1},\ldots,v_{j+l}), v_{j+l+1},
\ldots, v_n\right)=0.\nonumber\\ && \label{Phi_n}
\Eeqr

It is an easy calculation  using (\ref{la2}) to check that
\Beqrn
-\sum_{k+l=n+1\atop k,l\geq 2}
(-1)^{l+k(\tv_1+\ldots+\tv_l)}\la_l(v_1,\ldots,v_l)\cdot
[Q\la_{k-1}(v_{l+1},\ldots,v_n)] &&\\
+\sum_{k+l=n+1\atop k+l\geq 2} (-1)^{l(\tv_1+\ldots + \tv_{k-1})}
[Q\la_{k-1}(v_1,\ldots,v_{k-1})]\cdot \la_l(v_k,\ldots,v_n)&=&0.
\Eeqrn
Then  the first two sums in  (\ref{Phi_n})
reduce, again with the help of (\ref{la2}), to the following expression
$$
-\sum_{k+l=n+1\atop k,l\geq 2}
\sum_{{s+t=k\atop s\geq 2}\atop
t\geq 1} (-1)^p
\left[Q\la_s\left(\la_l(v_1,\ldots,v_l), v_{l+1},\ldots,v_{s+l-1}
)\right)\right]\cdot
\left[Q\la_t(v_{l+s},\ldots,v_n)\right]
$$
$$
-\sum_{k+l=n+1\atop k,l\geq 2}
\sum_{{s+t=k\atop s\geq 1}\atop
t\geq 2} (-1)^q
\left[Q\la_s\left(v_1,\ldots, v_{s}\right)\right]\cdot
\left[Q\la_t\left(v_{s+1},\ldots,v_{k-1},\la_l(v_k,\ldots,v_n)\right)\right],
$$
where $p={k(l-1)+s+(t-1)(\tv_1+\ldots +\tv_{s+l-1}+l)}$ and
$q=l(\tv_1+\ldots+\tv_{k-1})+k-1 + s + (t-1)(\tv_1+\ldots+\tv_s)$.
The third sum in (\ref{Phi_n}) splits into the following two sums
$$
\sum_{k+l=n+1\atop k,l\geq 2}\sum_{j=1}^{k-2}\sum_{{s+t=k\atop
1\leq s\leq j}\atop t\geq 1}(-1)^a
\left[Q\la_s(v_1,\ldots,v_s)\right]\cdot \left[Q\la_t\left(
v_{s+1},\ldots, v_j,\la_l(v_{j+1},\ldots,v_{j+l}), v_{j+l-1},\ldots,v_n\right)
\right]
$$
$$
-\sum_{k+l=n+1\atop k,l\geq 2}\sum_{j=1}^{k-2}\sum_{{s+t=k\atop
s\geq j+1}\atop t\geq 1}(-1)^b
\left[Q\la_s\left(
v_1,\ldots, v_j,\la_l(v_{j+1},\ldots,v_{j+l}), v_{j+l+1},\ldots,v_{s+l-1}
\right)
\right]\cdot \left[Q\la_t(v_{s+l},\ldots,v_n)\right],
$$
where $a=r + s-1 + (t-1)(\tv_1+\ldots+\tv_s)$ and
$b=r + s + (t-1)(\tv_1+\ldots+\tv_{s+l-1}+l)$.

\vst

Substituting these two expressions back into (\ref{Phi_n}) one easily gets
the following recursive formula
\Beqrn
\Phi_n(v_1,\ldots,v_n)&=& \sum_{{k+l=n \atop k\geq 3}\atop l\geq 1}
(-1)^{(l-1)(\tv_1+\ldots+\tv_k)+k}\left[Q\Phi_k(v_1,\ldots,v_k)\right]
\cdot \left[Q\la_l(v_{k+1},\ldots,v_n)\right]\\
&& - \sum_{{k+l=n \atop k\geq 1}\atop l\geq 3}
(-1)^{l(\tv_1+\ldots+\tv_k)}\left[Q\la_k(v_1,\ldots,v_k)\right]
\cdot \left[Q\Phi_l(v_{k+1},\ldots,v_n)\right],
\Eeqrn
where $n\geq 4$.

\vst

Finally we  compute
$$
\Phi_3(v_1,v_2,v_3)= (v_1\cdot v_2)\cdot v_3 - v_1\cdot (v_2\cdot v_3)=0.
$$
Thus $\Phi_n=0$ for all $n\geq 3$. $\Box$

\vsf

\begin{lemma}\label{lemma2}
The tensors $\la_k$, $k\geq 2$, satisfy,
for any $n\geq 2$ and any $v_1,\ldots,v_n\in V$, the identities
\Beqrn
\Theta_n(v_1,\ldots,v_n)&\equiv& d\la_n(v_1,\ldots,v_n) +
\sum_{j=0}^{n-1}(-1)^{n-1+\tv_1+\ldots+\tv_j}\la_n(v_1,\ldots, v_j,
dv_{j+1},v_{j+2},\ldots,v_n) \\
&&-
\sum_{k+l=n+1\atop k,l\geq 2}\sum_{j=0}^{k-1}
(-1)^r \la_k\left(v_1,\ldots,v_j, [d,Q]
\la_l(v_{j+1},\ldots,v_{j+l}), v_{j+l+1},
\ldots, v_n\right)\\
&=&0
\Eeqrn
where $r=l(\tv_1+\ldots \tv_j) + j(l-1) + (k-1)l$.

\end{lemma}
\noindent{\bf A scetch of the proof}. Following the same scenario
as the one used in the proof of Lemma~\ref{lemma1} (i.e.\ studing
separately the group of terms with extreme values of the index $j$, where
some cancellations occur, and applying (\ref{la2}) throughout), one gets
the following recursion formula
\Beqrn
\Theta_n(v_1,\ldots,v_n) &=&
 \sum_{{k+l=n \atop k\geq 2}\atop l\geq 1}
(-1)^{(l-1)(\tv_1+\ldots+\tv_k)+k}\left[Q\Theta_k(v_1,\ldots,v_k)\right]
\cdot \left[Q\la_l(v_{k+1},\ldots,v_n)\right]\\
&& - \sum_{{k+l=n \atop k\geq 1}\atop l\geq 2}
(-1)^{l(\tv_1+\ldots+\tv_k)}\left[Q\la_k(v_1,\ldots,v_k)\right]
\cdot \left[Q\Theta_l(v_{k+1},\ldots,v_n)\right],
\Eeqrn
where $n\geq 3$.

\vst

Next we compute
$$
\Theta_2(v_1,v_2)= d(v_1\cdot v_2) - (dv_1)\cdot v_2 - (-1)^{\tv_1}
v_1\cdot (dv_2)= 0.
$$
Thus $\Theta_n=0$ for all $n\geq 2$. $\Box$

\vsf

\vst

\begin{theorem}\label{main}
Let $(V,d)$ be a differential associative superalgebra and let $(W,d)\subset
(V,d)$ be a subcomplex satisfying the Assumption~\ref{assumption}.
Then the linear maps
$$
\mu_k: \ot^k W \lon W, \ \ \  k\geq1,
$$
defined by
\Beqr
\mu_1 &:=& d,\nonumber\\
\mu_k &:=& (1-[d,Q])\la_k, \ \ \mbox{for}\ {k\geq 2}, \label{mu}
\Eeqr
with $\la_k$ being given by (\ref{la2}), satisfy the higher order
associativity identitites (\ref{id}) for all $n\geq 1$.

Thus there exists on  $W$
a structure of $A_{\infty}$-algebra.
\end{theorem}
\noindent{\bf Proof}. Denote the l.h.s. of the equation (\ref{id}) by
$\Psi_n$. Then $\Psi_1$ and $\Psi_2$ vanish because $(W,d)$ is a subcomplex
of $(V,d)$. For $n\geq 3$ one has
$$
\Psi_n = (1-[d,Q])\left(\Phi_n + \Theta_n\right).
$$
By Lemmas~\ref{lemma1} and \ref{lemma2}, the tensors $\Phi_n$ and
$\Theta_n$ vanish for all $n\geq 3$. This completes the proof. $\Box$

\vst

The existence part of Theorem~\ref{main} has been proved  by
Gugenheim and Stasheff \cite{GS} under slightly stronger
assumption and using different method. They built up an inductive series of
approximations $\mu^{(n)}_k$, $n\in {\Bbb N}$, to the higher
homotopies $\mu_k$ and proved its
convergence in the limit $n\rar \infty$
to a genuine $A_{\infty}$-structure\footnote{After this work was completed
L.\ Johansson drew author's attention to the paper \cite{JL} where the
limits
$\lim_{n\rar\infty}\mu_k^{(n)}$ have been  computed.}.

\vsf

Finally we note that the higher homotopies $\mu_n$, $n\geq 2$,
can be written as follows
$$%\label{la2}
\mu_n(v_1,\ldots,v_n) = -\sum_{ k+l=n+1\atop k,l\geq 1 }(-1)^{k+(l-1)(\tv_1+\ldots+
\tv_k)} [Q\la_k(v_1, \ldots, v_k)] \circ [Q\la_l(v_{k+1},\ldots,v_n)].
$$

%%%%%%%%%%%%%%%%%%%%%%%%%%%%%%%%%%%%%%%%%%%%%%%%%%%%%%%%
\section{ $A_{\infty}$-algebras of a K\"ahler manifold}

\paragraph{\bf 4.1. Hodge theory: a brief overview.}
Let $M$ be a compact
K\"ahler manifold, $\Omega^{\bullet}M$ the algebra of real smooth
 differential forms on $M$, $\Omega^{\bullet,\bullet}M$ the algebra
of complex differential forms with Hodge gradation, and
$$
\p:\Omega^{\bullet,\bullet}M \lon  \Omega^{\bullet+1,\bullet}M, \ \ \ \ \ \
\bp:\Omega^{\bullet,\bullet}M \lon  \Omega^{\bullet,\bullet+1}M,
$$
the standard derivations. We also consider two real derivations
$$
d=\p+\bp: \Omega^{\bullet}M \lon  \Omega^{\bullet+1}M, \ \ \ \ \
d_c=i(\p-\bp):\Omega^{\bullet}M \lon  \Omega^{\bullet+1}M.
$$

\vst

The K\"ahler metric $g$ on $M$
gives rise to a Hermitian metric on $\Omega^{\bullet,\bullet}$ defined,
say for $(p,q)$-forms
$$
\al= \sum\al_{i_1\ldots i_p\bbj_1\ldots \bbj_q}dz^{i_1}\wedge \ldots
\wedge dz^{i_p} d\bbz^{\bbj_1}\wedge \ldots
\wedge d\bbz^{\bbj_q},
$$
and
$$
\be= \sum\be_{s_1\ldots s_p\bbt_1\ldots \bbt_q}dz^{s_1}\wedge \ldots
\wedge dz^{s_p} d\bbz^{\bbt_1}\wedge \ldots
\wedge d\bbz^{\bbt_q},
$$
by the integral
$$
(\al,\be):=\int_M \left(\sum g^{i_1\bbs_1}\ldots g^{i_p\bbs_p}
g^{t_1\bbj_1}\ldots g^{t_q\bbj_q}\,
\al_{i_1\ldots i_p\bbj_1\ldots \bbj_q}\,
\overline{\be_{s_1\ldots s_p\bbt_1\ldots \bbt_q}}\right) \mbox{vol}_g,
$$
where $g_{i\bbj}$ are components of $g$ in a local holomorphic coordinate
chart $\{z^i\}$. The scalar product $(\al,\be)$ of differential
forms with different Hodge type is set to be zero.

\vst

This scalar product is used to define conjugate operators
$$
\p^*:\Omega^{\bullet,\bullet}M \lon  \Omega^{\bullet-1,\bullet}M,  \ \ \
\bp^*:\Omega^{\bullet,\bullet}M \lon  \Omega^{\bullet,\bullet-1}M, \ \
\Lambda: \Omega^{\bullet,\bullet}M \lon  \Omega^{\bullet-1,\bullet-1}M, \ \
$$
by the formulae
$$
(\p\al, \be)=(\al, \p^* \al), \ \ \ (\bp\al,\be)=(\al,\bp^*\be),
\ \ \ (\al\wedge \omega,\be)=(\al,\Lambda\be),
$$
where $\omega$ denotes the K\"ahler form on $M$ and $\al,\be\in
\Omega^{\bullet,\bullet}M$ are arbitrary. We shall also need
the following real operators
$$
d^* = \p^* + \bp^*, \ \ \ \ \ d_c^*=-i(\p^*-\bp^*).
$$

There are four Laplacians,
$$
\Delta_{\p}= \p\p^*+\p^*\p, \ \ \ \ \
\Delta_{\bp}= \bp\bp^*+\bp^*\bp,
$$
$$
\Delta_{d}= dd^*+d^*d, \ \ \ \ \
\Delta_{d_c}= d_c d_c^*+ d_c^*d_c.
$$
One of the central results of the Hodge theory says that
$$
\Delta_{\p}=\Delta_{\bp}=\frac{1}{2}\Delta_d=\frac{1}{2}\Delta_{d_c}.
$$

A differential form $\al\in \Omega^{\bullet,\bullet}M$ is called {\em harmonic}
if $\Delta_d \al=0$. The vector space Harm of all real (or complex)
harmonic forms on $M$ is finite dimensional and is isomorphic to
its de Rham (or Dolbeault) cohomology. There is a natural
orthogonal projection
$$
\Ba{rccc}
[\ ]_{\mathrm Harm}: & \Omega^{\bullet,\bullet}M & \lon & \mbox{Harm}\\
& \al & \lon & [\al]_{\mathrm Harm}.
\Ea
$$
Note that all the Laplacians $\Delta$ become isomorphisms when restricted
to the orthogonal
complement, $\mbox{Harm}^{\perp}$, of the space of harmonic forms.

\vst

There are very important for our purposes Hodge identities
\Beq\label{rel}
[\Lambda,\bp]=-i\p^*, \ \ \ \ \  [\Lambda,\bar{\p}]=i\p^*, \ \ \ \ \
[\Lambda,d]=d_c^*,
\Eeq
and the Hodge decompositions, for an arbitrary
$\al\in \Omega^{\bullet,\bullet}M$,
\Beqrn
\al&=&[\al]_{\mathrm Harm} + \p G_{\p} \p^*(\al) + \p^* G_{\p} \p(\al)\\
\al&=&[\al]_{\mathrm Harm} + \bp G_{\bp} \bp^*(\al) + \bp^* G_{\bp} \bp(\al)\\
\al&=&[\al]_{\mathrm Harm} + d G_{d} d^*(\al) + d^* G_{d} d(\al)\\
\al&=&[\al]_{\mathrm Harm} + d_c G_{d_c} d_c^*(\al) + d_c^* G_{d_c} d_c(\al),
\Eeqrn
where $G: \Omega^{\bullet,\bullet}M\lon \Omega^{\bullet,\bullet}M$
are the Green operators defined by
$$
G_{\p}\mid_{\mathrm Harm}=0, \ \ \ \ \  G_{\p}\mid_{{\mathrm Harm}^{\perp}}=
\Delta_{\p}^{-1},
$$
and analogously for all others. A classical fact:
$$
G_{\p}=G_{\bp}= 2G_{d}=2G_{d_c}.
$$
Note that $\p G_{\p}= G_{\p}\p$,
$\p^* G_{\p}= G_{\p}\p^*$, etc.

\vst

For further details about Hodge theory on K\"ahler manifolds we refer to
\cite{GH}.

\vst

\paragraph{\bf 4.2. Real $A_{\infty}$-algebra of a K\"ahler manifold.}
Consider the subcomplex $(W:=\Ker d_c^*, d)$ in the de Rham differential
superalgebra $(\Omega^{\bullet}M,d)$. It does satisfy the
Assumption~\ref{assumption} with
$$
Q=d_c G_d \Lambda.
$$
Indeed,
$$
1-[d,Q]= 1- d_c G_d [\Lambda,d]= 1-d_c G_d d_c^*
$$
which, according to the Hodge decomposition, is precisely
the projector from $\Omega^{\bullet}M$ to its subspace $\Ker d_c^*$.

\vst

By Theorem~\ref{main}, the space $\Ker d_c^*$ has canonically
the structure of an $A_{\infty}$-algebra with the higher homotopies
given explicitly by the formulae (\ref{mu}) and (\ref{la1}). For example,
\Beqrn
\mu_1 &=& d,\\
\mu_2(v_1,v_1) &=& (1-d_cG_d d_c^*)(v_1\wedge v_2)\equiv v_1\circ v_2,\\
\mu_3(v_1,v_2,v_3) &=& \left[d_cG_d \Lambda(v_1\wedge v_2)\right]\circ v_3
- (-1)^{\tv_1}v_1\circ \left[d_c G_d \Lambda (v_2\wedge v_3)\right], \\
\mbox{etc.} &&
\Eeqrn

In this algebra the product of two harmonic forms is obviously
harmonic as well. Another interesting observation is that
the complex $(\Ker d_c^*, d)$ is precisely the one which gives rise
to a Frobenius manifold structure on the de Rham cohomology
of $M$ \cite{Me}. The same comments will apply to our next example.

\vst

\paragraph{\bf 4.3. Complex $A_{\infty}$-algebra of a K\"ahler manifold.}
Consider the subcomplex $(W:=\Ker \p^*, \bp)$ in the Dolbeault
differential superalgebra $(\Omega^{\bullet,\bullet}M,\bp)$.
It satisfies the
Assumption~\ref{assumption} with
$$
Q=i\p G_{\p}\Lambda.
$$
Indeed,
$$
1-[d,Q]= 1- i\p G_{\p} [\Lambda,\bar{\p}]= 1-\p G \p^*
$$
which, according to the Hodge decomposition, is precisely
the projector from $\Omega^{\bullet,\bullet}M$ to the
subspace $\Ker \p^*$.

\vst

Therefore, by Theorem~\ref{main}, the space $\Ker \p^*$ has canonically
the structure of an $A_{\infty}$-algebra
and one may explicitly write all the higher homotopies $\mu_k$
in terms of operators $\p$, $\p^*$, $G_{\p}$ and $\Lambda$.

\vst

\paragraph{\bf 4.4. $A_{\infty}$-algebra of a Calabi-Yau manifold.}
A compact K\"ahler manifold $M$ is called Calabi-Yau if $c_1(M)=0$
or, equivalently, if it admits a nowhere vanishing holomorphic
section $\Omega$ of the canonical line bundle of $M$.

\vst

Consider a differential associative superalgebra
$$
V:= \Lambda^{\bullet} TM\ot \Lambda^{\bullet}\overline{\Omega^1 M}
%=\sum_{p,q} \Gamma(M, \Lambda^{p} TM
%\ot \Lambda^q \overline{\Omega^1 M},
$$
with the differential $\bp$ and the multiplication \, $\cdot$ \,
being given by the ordinary wedge products on $\Lambda^{p} TM$ and
$\Lambda^q \overline{\Omega^1 M}$. It was first considered
by Barannikov and Konstevich \cite{BaKo} in the context of
mirror symmetry.

\vst

The holomorphic volume form
defines an isomorphism
$$
\Ba{rccc}
\Omega:& \Lambda^{\bullet} TM \ot\Lambda^{\bullet} \overline{\Omega^1 M}
& \lon & \Lambda^{\dim M -\bullet} \Omega^1 M \ot
\Lambda^{\bullet} \overline{\Omega^1 M}\\
& \al\ot\be& \lon & (\al\lrcorner\, \Omega)\ot \be.
\Ea
$$
Then for any linear operator
$$
s:\Omega^{\bullet,\bullet}M\lon
\Omega^{\bullet,\bullet}M
$$
there is associated a linear operator
$$
\hat{s}: V\lon V
$$
defined as the composition:
$$
\hat{s}: V\stackrel{\Omega}{\lon} \Omega^{\bullet,\bullet}M
\stackrel{s}{\lon} \Omega^{\bullet,\bullet}M
\stackrel{\Omega^{-1}}{\lon} V.
$$
In this way one obtains the operators $\hat{\p}$, $\hat{\p}^*$,
$\hat{G}_{\p}$, etc. Note that $\hat{\bp}=\bp$. Note also that
$\hat{\p}$ is not a derivation of $(V,\cdot)$.

\vst

Consider the following subcomplex
$$
(W:=\Ker \hat{\p}, \bp)
$$
 in the Barannikov-Kontsevich differential
superalgebra $(V,\bp)$. It does satisfy the
Assumption~\ref{assumption} with
$$
Q=-i\hat{\Lambda} \hat{G}_{\p}\hat{\p}.
$$
Indeed,
$$
1-[d,Q]= 1- i\widehat{[\Lambda,\bp]}\hat{G}_{\p}\hat{\p} =
 1-\hat{\p}^*\hat{G}_{\p}\hat{\p},
$$
which, according to the Hodge decomposition, is precisely
the projector from $V$ to the
subspace $\Ker \hat{\p}$.

\vst

Therefore, by Theorem~\ref{main}, the space $\Ker \hat{\p}$ has canonically
the structure of an $A_{\infty}$-algebra
and one may explicitly write down all the higher homotopies $\mu_k$
as explained in Sect.3.
It is again a remarkable fact that the same complex $(\Ker \hat{\p},\bp)$
admits not only an $A_{\infty}$ but also a very special
Gerstenhaber-Batalin-Vilkoviski structure which makes the extended moduli space
of complex structures on $M$ into a Frobenius manifold \cite{BaKo}.

\vse

{\small {\em Acknowledgments}. This work has been done during author's visit
to the Max-Planck Institute for Mathematics in Bonn. It is a pleasure
to thank Yu.I.\ Manin, V.\ Schechtman and J.\ Stasheff for
valuable discussions and comments.}

\vse

\pagebreak
{\small

{\sc
\begin{tabular}{l}
Department of Mathematics\\
University of Glasgow\\
15 University Gardens \\
 Glasgow G12 8QW, UK\\
 \mbox{\rm sm@maths.gla.ac.uk}
\end{tabular}
}

\end{document}